\documentclass[11pt]{article}

\usepackage{a4,amssymb,amsmath,color,graphicx}
\usepackage{trfsigns}

\def\eins{\mbox{1\hskip-0.24em l}}
\def\T{^{\sf T}}

\newcommand{\N}{ {\mathbb N} }
\newcommand{\R}{ {\mathbb R} }

\newcommand{\diag}{\,\mbox{diag}}

\newcommand{\qed}{\qquad\mbox{$\square$}}

\newtheorem{remark}{Remark}[section]

\newtheorem{lemma}{Lemma}[section]

\begin{document}
\title{A Stiff MOL Boundary Control Problem for the 1D Heat Equation
with Exact Discrete Solution}
\author{Jens Lang \\
{\small \it Technical University of Darmstadt,
Department of Mathematics} \\
{\small \it Dolivostra{\ss}e 15, 64293 Darmstadt, Germany}\\
{\small lang@mathematik.tu-darmstadt.de} \\ \\
Bernhard A. Schmitt \\
{\small \it Philipps-Universit\"at Marburg,
Department of Mathematics,}\\
{\small \it Hans-Meerwein-Stra{\ss}e 6, 35043 Marburg, Germany} \\
{\small schmitt@mathematik.uni-marburg.de}}
\maketitle

\begin{abstract}
Method-of-lines discretizations are demanding test problems for stiff integration methods.
However, for PDE problems with known analytic solution the presence of space discretization errors or the need to use codes to compute reference solutions may limit the validity of numerical test results.
To overcome these drawbacks we present in this short note a simple test problem with boundary control, a situation where one-step methods may suffer from order reduction. We derive exact formulas for the
solution of an optimal boundary control problem governed by a one-dimensional discrete heat equation and an objective function that measures the distance of the final state from the target and the control costs. This
analytical setting is used to compare the numerically observed convergence orders for selected implicit Runge-Kutta and Peer two-step methods of classical order four which are suitable for optimal control problems.
\end{abstract}

\noindent{\em Key words.} Optimal control, implicit Peer two-step methods, first-discretize-then-optimize, discrete adjoints

\section{Introduction}
One main area of application for stiff integration methods is semi-discretizations in space of time-dependent partial differential equations in the method-of-lines approach.
In order to test new methods in this area one may rely on PDE test problems with known analytic solution or reference solutions computed by other numerical methods.
However, both approaches have its drawbacks.
For PDE problems the accuracy is limited by the level of space discretization errors and in computing reference solutions one has to trust the reliability of the used code.
This background was our motivation to develop the current test problems with exact discrete solutions for a finite difference semi-discretization in space with arbitrarily fine grids.
\par
It is known that, in contrast to multi-step type methods, one-step methods may suffer from order reduction if applied to MOL systems especially with time-dependent boundary conditions, see \cite{LubichOstermann1995,OstermannRoche1992}.
Motivated by our recent work \cite{LangSchmitt2022a,LangSchmitt2022b} on Peer two-step methods in optimal control the present example is formulated as a problem with boundary control.
\par
The paper is organized as follows. In Section 2, we apply a finite difference discretization with a shifted
equi-spaced grid for the 1D heat equation with general Robin boundary conditions and derive exact formulas
for the solutions of the discrete heat equation and an optimal boundary control problem. These analytical
solutions are used in a sparse setting to study the numerically observed convergence orders in Section 3 for several one-step and two-step integration methods which are suitable for optimal control.
\section{A discrete heat equation with boundary control}
\subsection{Finite difference discretization of the 1D heat equation}
We consider the initial-boundary-value problem for a function $Y(x,t)$ governed by the heat equation
\begin{align}
\label{WLGl}
\partial_tY(x,t) = &\,\partial_{xx}Y(x,t),\ (x,t)\in[0,1]\times[0,T],\\[1mm]
\label{RRBu}
\partial_xY(0,t) = &\,0,\; \beta_0 Y(1,t)+\beta_1 \partial_xY(1,t)=u(t),\\[1mm]
\label{ABPDG}
Y(x,0)= &\, \Psi(x),
\end{align}
where $\Psi(x)$ and $u(t)$ are given functions.
The homogeneous Neumann condition at $x=0$ may be considered as a short-cut
for space-symmetric solutions $Y(-x,t)\equiv Y(x,t)$.
The coefficients of the general Robin boundary condition are nonnegative,
$\beta_0,\beta_1\ge0$ and nontrivial $(\beta_0,\beta_1)\not=(0,0)$.
\par
The equation \eqref{WLGl} is approximated by finite differences with a shifted
equi-spaced grid with stepsize $\xi=1/m$, $m\in\N$:
\begin{align}\label{grid}
 x_j=\left( j-\frac12\right)\xi,\ j=1,\ldots, m.
\end{align}
For the approximation of the boundary conditions, also the outside points $x_0=-\xi/2$
and $x_{m+1}=1+\xi/2$ will be considered temporarily.
In the method-of-lines approach with central differences, approximations
$y_j(t),\,j=1,\ldots,m$ are defined by the differential equations
\begin{align}\label{MOLj}
 y_j'=\frac1{\xi^2}\left( y_{j-1}-2y_j+y_{j+1}\right),\ j=2,\ldots,m-1,
\end{align}
for the grid points in a distance to the boundary.
The symmetric difference approximation $0\stackrel!=\xi Y_x(0,t)\cong(y_1-y_0)$ leads
to the symmetry condition $y_0\equiv y_1$ and yields the MOL-equation
\begin{align}\label{MOL1}
y_1'=\frac{-y_1+y_2}{\xi^2}.
\end{align}
In a similar way, the Robin boundary condition is approximated by the equation
\begin{align}
\beta_0\frac{y_m+y_{m+1}}2+\beta_1\frac{y_{m+1}-y_m}{\xi}=u(t),
\end{align}
which may be solved for $y_{m+1}$ by
\begin{align}
 y_{m+1}=\frac{2\beta_1-\beta_0\xi}{2\beta_1+\beta_0\xi}y_m+
 \frac{2\xi}{2\beta_1+\beta_0\xi}u(t).
\end{align}
Thus, $y_{m+1}$ may be eliminated from the equation \eqref{MOLj} with $j\!=\!m$ yielding
\begin{align}\label{MOLm}
 y_m'=&\frac1{\xi^2}(y_{m-1}-\theta y_m)+\gamma\,u(t)
\end{align}
with
\begin{align}\label{theta}
 \theta=\frac{2\beta_1+3\beta_0\xi}{2\beta_1+\beta_0\xi}=3-\frac{4\beta_1}{2\beta_1+\beta_0\xi},
 \quad \gamma=\frac{2}{(2\beta_1+\beta_0\xi)\xi}.
\end{align}
Hence, we have $\theta\!=\!3$ for the Dirichlet condition and $\theta\!=\!1$ for the pure Neumann condition.
Collecting all equations \eqref{MOLj}, \eqref{MOL1} and \eqref{MOLm}, the following MOL system for the vector $y(t)=\big(y_j(t)\big)_{j=1,\ldots,m}$ is obtained:
\begin{align}\label{MOLg}
 y'=&My+\gamma e_mu(t),\\[3mm]\label{MatM}
 M =&\frac1{\xi^2}\begin{pmatrix}
   -1&1\\
   1&-2&1\\
   &&\ddots&\ddots&\ddots\\
   &&&1&-2&1\\
   &&&&1&-\theta
 \end{pmatrix},
\end{align}
where $e_m$ is the $m$-th unit vector.
The initial conditions are simple evaluations of the function $\Psi$ on the grid,
\begin{align}\label{MOLAB}
 y(0)=\psi,\quad \psi=\left( \Psi(x_j)\right)_{j=1}^m.
\end{align}
\par
The basis of our construction is that the eigenvalues and eigenvectors of the
symmetric matrix $M$ are known, which is well known for special values of
$\theta$, at least.
\begin{lemma}\label{LEWV}
For $m\ge2$, the eigenvalues of the matrix $M\in\R^{m\times m}$ from \eqref{MatM}
are given by
\begin{align}
\lambda_k=-4m^2\sin^2\left( \frac{\omega_k}{2m}\right),\ k=1,\ldots,m,
\end{align}
where $\omega_k,\,k=1,\ldots,m,$ are the $m$ first non-negative solutions of the equation
\begin{align}\label{EWGl}
\tan(\omega)\tan\left( \frac{\omega}{2m}\right) =\frac{\beta_0}{2m\beta_1},
\end{align}
with the convention that $\omega_k=(k-\frac12)\pi$, $k=1,\ldots,m,$ for $\beta_1=0$.
The corresponding normalized eigenvectors $v^{[k]}$ have the components
\begin{align}\label{EVek}
 v_j^{[k]}=&\nu_k\cos\left( \omega_k\frac{2j-1}{2m}\right),\ j=1,\ldots,m,
\end{align}
with constants $\nu_k=2/\sqrt{2m+\sin(2\omega_k)/\sin(\omega_k/m)}$.
\end{lemma}
{\bf Proof:}
In the main equation \eqref{MOLj}, the ansatz $v=\big(\Re e^{i\omega x_j}\big)_{j=1}^m$
gives
\begin{align*}
 \frac1{\xi^2}(v_{j-1}-2v_j+v_{j+1})
 =&\frac1{\xi^2}\Re e^{i\omega x_j}\left( e^{-i\omega\xi)}-2+e^{i\omega\xi}\right)
  =-\frac4{\xi^2}\sin^2\left( \frac{\omega\xi}2\right) v_j.
\end{align*}
In the first equation, we have
\begin{align*}
\frac1{\xi^2}(-v_1+v_2)=&\,\frac1{\xi^2}\left(-\cos{\frac{\omega\xi}2}+
\cos\left(3\frac{\omega\xi}2\right)\right)
 =\frac4{\xi^2}\left(\cos^3\left(\frac{\omega\xi}2\right)-\cos{\frac{\omega\xi}2}\right)\\
 &\,=-\frac4{\xi^2}\sin^2\left(\frac{\omega\xi}2\right)v_1,
\end{align*}
with the same factor $\lambda\!:=\!-(4/\xi^2)\sin^2\big({\omega\xi}/2\big)$.
In order to satisfy the eigenvalue condition in the last component, we consider the equation $0=e_m\T(Mv-\lambda v)$, i.e.,
\begin{align*}
 0\stackrel!=&\,v_{m-1}-\left( \theta+\lambda\xi^2\right) v_m
  =\cos\left( \omega(x_m-\xi)\right)-\left( \theta+\lambda\xi^2\right)\cos(\omega x_m)\\[2mm]
  =&\,\left( \cos(\omega\xi)+4\sin^2\left( \frac{\omega\xi}2\right)-\theta\right) \cos(\omega x_m)+\sin(\omega\xi)\sin(\omega x_m)\\[2mm]
  =&\,\left(1 +2\sin^2\left( \frac{\omega\xi}2\right)-\theta\right)\cos(\omega x_m)+\sin(\omega\xi)\sin(\omega x_m),
\end{align*}
since $\cos(\omega\xi)\!=\!1-2\sin^2(\omega\xi/2)$.
The last grid point is $x_m\!=\!1-\xi/2$ and with the trigonometric formulas for $\cos(\omega-\omega\xi/2),\sin(\omega-\omega\xi/2)$ and the identity $\sin(\omega\xi)=2\sin(\omega\xi/2)\cos(\omega\xi/2)$, we may proceed with
\begin{align*}
0=&\,\left( 1+2\sin^2\left(\frac{\omega\xi}2\right)-\theta\right)
\left(\cos(\omega)\cos\left(\frac{\omega\xi}2\right)+
\sin(\omega)\sin\left(\frac{\omega\xi}2\right)\right)\\[2mm]
&\,+2\sin\left(\frac{\omega\xi}2\right)\cos\left(\frac{\omega\xi}2\right)
\left(\sin(\omega)\cos\left(\frac{\omega\xi}2\right)
-\cos(\omega)\sin\left(\frac{\omega\xi}2\right)\right)\\[2mm]
=&\,(1-\theta)\cos\left(\frac{\omega\xi}2\right)\cos(\omega)
+\left( 1+2\sin^2\left(\frac{\omega\xi}2\right)-\theta+2\cos^2\left(\frac{\omega\xi}2\right)\right)
\sin\left(\frac{\omega\xi}2\right)\sin(\omega)\\[2mm]
=&\,(1-\theta)\cos\left(\frac{\omega\xi}2\right)\cos(\omega)
+(3-\theta)\sin\left(\frac{\omega\xi}2\right)\sin(\omega).
\end{align*}
Hence, the different versions of $\theta$ in \eqref{theta} verify the condition \eqref{EWGl}
for $m=1/\xi$.
Rearranging \eqref{EWGl} as $\tan(\omega)=\beta_0/(2m\beta_1)\cot(\omega/(2m))$, for $\beta_0>0$ it is seen that exactly $m$ solutions exist in $(0,m\pi)$ since the function $\omega\mapsto\cot(\omega/(2m))$ is monotonically decreasing and positive.
Finally, the vector norms are computed for $\omega\not=0$.
Abbreviating $\omega/m=:\Omega$ and using $\cos^2(x)=(1+\cos(2x))/2$, we get
\begin{align*}
\sum_{j=1}^m\cos^2\left(\left(j-\frac12\right)\frac\omega m\right)
=&\,\frac{m}2+\frac12\sum_{j=1}^m\cos\left( (2j-1)\Omega\right)\\[1mm]
=&\,\frac{m}2+\frac12\Re\sum_{j=1}^me^{i(2j-1)\Omega}
=\frac{m}2+\frac12\Re\frac{e^{i2\Omega m}-1}{e^{i\Omega}-e^{-i\Omega}}\\[1mm]
=&\,\frac{m}2+\frac14\Im\frac{e^{i2\Omega m}-1}{\sin(\Omega)}
=\frac{m}2+\frac14\frac{\sin(2\omega)}{\sin(\Omega)},
\end{align*}
which leads to the value of the normalizing factor $\nu$ in \eqref{EVek}.
\qed
\par
For later use, we introduce the diagonal matrix $\Lambda=\diag(\lambda_k)$ and
the unitary matrix $V=(v^{[1]},\ldots,v^{[m]})$ satisfying $M=V\Lambda V\T$.
\par
\begin{remark}
The well-known frequencies for Dirichlet boundary conditions are $\omega_k=(k-\frac12)\pi$
and  $\omega_k=(k-1)\pi$, $k=1,\ldots,m$ for Neumann conditions.
For general values $\beta_0,\beta_1>0$ the equation \eqref{EWGl} may be rewritten in fixed point form
\begin{align}
\omega=f_k(\omega):=(k-1)\pi+\arctan\left(\frac{\beta_0}{2m\beta_1}\cot\left(\frac{\omega}{2m}\right)\right),
\quad k=1,\ldots,m.
\end{align}
The functions $f_k$ are monotonically decreasing in $\omega$ and an iteration with initial value $\omega=\max\{1,(k-1)\pi\}$ converges to the desired solution $\omega_k$ at least for $2m>\beta_0/\beta_1\ge1$, since $1/|f_k'(\omega)|=(4m^2\beta_1/\beta_0)\sin^2(\omega/(2m))+(\beta_0/\beta_1)\cos^2(\omega/(2m))$.
\end{remark}
\subsection{Exact solution of the discrete heat equation}
Knowing the eigenvectors and eigenvalues of the linear problem \eqref{MOLg}, the computation of its solution is straight-forward.
The representation $y(t)=\sum_{k=1,\ldots,m}\eta_k(t)v^{[k]}$ leads to
\begin{align}
\sum_{k=1}^m\eta_k'(t)v^{[k]}=\sum_{k=1}^m\lambda_k\eta_k(t) v^{[k]}+\gamma e_mu(t).
\end{align}
Since the matrix $M$ is symmetric, the inner product with $v^{[j]}$ yields the decoupled equations $\eta_j'(t)=\lambda_j\eta_j(t)+\gamma v_m^{[j]}u(t)$, which can be solved easily leading to the following result.
\begin{lemma}\label{LMOLsg}
With the data from Lemma~\ref{LEWV}, the solution of the initial value problem \eqref{MOLg},
\eqref{MOLAB} is given by
\begin{align}\label{LsgAWP}
 y(t)=&\,\sum_{k=1}^m\left( e^{\lambda_kt}{v^{[k]}}\T \psi
  +\gamma v_m^{[k]}\int_0^t e^{\lambda_k(t-\tau)}u(\tau)\,d\tau\right)v^{[k]}.
\end{align}
\end{lemma}
\begin{remark} The presence of the terms $v_m^{[k]}$ indicates that simple sparse solutions with only a few terms in \eqref{LsgAWP} may not exist due to the inhomogeneous boundary condition \eqref{RRBu}.
\end{remark}
\subsection{Exact solution of an optimal control problem}
The inhomogeneity $u(t)$ inherited from the boundary condition \eqref{RRBu} may be considered as a control to approach a given target profile $\hat y\in\R^m$ at some given time $T>0$.
In an optimal control context, controls are searched for minimizing an objective function like
\begin{align}\label{Zielf}
 C=\frac12\|y(T)-\hat y\|_2^2+\frac{\alpha}2\int_0^T u(t)^2dt,
\end{align}
with $\alpha>0$.
Optima may be computed by using some multiplier function $p(t)$ for the ODE
restriction \eqref{MOLg} and considering the Lagrangian
\begin{align}\notag
L:=&\,-C+\int_0^Tp\T\bigl(y'-My-\gamma e_mu\bigr)\,dt+p\T(0)(y(0)-\psi) \\\label{Lagr}
=&\,-C
-\int_0^T\left( \left( p'\right)\T y+p\T\bigl( My+\gamma e_mu \bigr)\right)dt+p\T(T)y(T)-p\T(0)\psi.
\end{align}
The partial derivatives of the Lagrangian $L$ with respect to $p(t)$ and $p(T)$ recover
\eqref{MOLg}, \eqref{MOLAB} and the other ones are
\begin{align}
\label{KKTyt}
 \partial_{y(t)}L=&\,-p'-M\T p=-p'-M p,\\
 \label{KKTyT}
 \partial_{y(T)}L=&\,\hat y-y(T)+p(T)\\
 \label{KKTu}
 \partial_{u(t)}L=&\,-\alpha u-\gamma e_m\T p.
\end{align}
Hence, the Karush-Kuhn-Tucker conditions, $\partial_{(\cdot)}L=0$ in \eqref{KKTyt}--\eqref{KKTu},
show that the control $u(t)$ may be eliminated by
\begin{align}\label{uausp}
 u(t)=-\frac{\gamma}{\alpha} e_m\T p(t)=-\frac{\gamma}{\alpha}p_m(t),
\end{align}
and a necessary condition for the optimal solution is that it solves the
following boundary value problem:
\begin{align}\label{RWPv}
 y'=&\,My-\frac{\gamma}{\alpha}e_me_m\T p,\quad y(0)=\psi,\\\label{RWPa}
 p'=&\,-Mp,\quad p(T)=y(T)-\hat y.
\end{align}
The homogeneous differential equation \eqref{RWPa} for $p$ has the simple solution
\begin{align}\label{Lsgp}
 p(t)=e^{(T-t)M}p(T)=\sum_{\ell=1}^m e^{\lambda_\ell(T-t)}v^{[\ell]}{v^{[\ell]}}\T(y(T)-\hat y)
\end{align}
with the matrix exponential
\begin{align}
e^{sM} = &\,V \diag\left( e^{\lambda_ks}\right)V\T, \quad 0\le s\le T.
\end{align}
With $u$ given by \eqref{uausp}, this solution may be used in \eqref{LsgAWP}
to yield the solution for \eqref{RWPv} with the coefficient functions
\begin{align}\notag
\eta_k(t)=&\,e^{\lambda_k t}\eta_k(0)-\frac{\gamma^2}{\alpha}v_m^{[k]}\int_0^t
e^{\lambda_k(t-\tau)}p_m(\tau)d\tau
\\[1mm]\label{LsgOC}
=&\,e^{\lambda_k t}\eta_k(0)-\frac{\gamma^2}{\alpha}v_m^{[k]}\sum_{\ell=1}^m
v_m^{[\ell]}\int_0^te^{\lambda_k t+\lambda_\ell T-(\lambda_k+
\lambda_\ell)\tau}d\tau\cdot {v^{[\ell]}}\T(y(T)-\hat y).
\end{align}
Considering the function $\varphi_1(z)=\int_0^1 e^{zt}dt$ satisfying
$\varphi_1(z)=(e^z-1)/z$ for $z\not=0$ and $\varphi_1(0)=1$, the integral may
be written as
\begin{align}
\int_0^te^{\lambda_k t+\lambda_\ell T-(\lambda_k+\lambda_\ell)\tau}d\tau
=te^{\lambda_\ell(T-t)}\varphi_1\left( (\lambda_k+\lambda_\ell)t\right).
\end{align}
The result \eqref{LsgOC} may be used in different ways.
The first application is computing the solution for a given target profile $\hat y$.
\begin{lemma}\label{LLsgtp}
Let $\hat y\in\R^m$ be given.
Then, the coefficient vector $\eta(T)$ of the solution $y(t)=V\eta(t)$ of the
boundary value problem \eqref{RWPv}, \eqref{RWPa} is given by the unique solution
of the linear system
\begin{align}\label{LGSend}
 (I+Q)\eta(T)=e^{T\Lambda}\eta(0)+QV\T\hat y,
\end{align}
with the positive semi-definite matrix $Q=(q_{k\ell})_{k,\ell=1}^m$ having the elements
\begin{align}\label{MatQ}
 q_{k\ell}=\frac{\gamma^2T}{\alpha}v_m^{[k]}\varphi_1\left(
 (\lambda_k+\lambda_\ell)T\right)v_m^{[\ell]},\ k,
\ell=1,\ldots,m.
\end{align}
\end{lemma}
{\bf Proof:}
At the end point $T$, the formula \eqref{LsgOC} simplifies to
\begin{align*}
\eta_k(T)=&\,e^{\lambda_k T}\eta_k(0)-\frac{\gamma^2T}{\alpha}v_m^{[k]}\sum_{\ell=1}^m
v_m^{[\ell]}\varphi_1\left( (\lambda_k+\lambda_\ell)T\right)\cdot
\left(\eta_\ell(T)-{v^{[\ell]}}\T\hat y\right).
\end{align*}
This equation may be reordered to the form given in \eqref{LGSend} with the matrix elements \eqref{MatQ}.
Finally, we consider the quadratic form of the matrix $Q$ with some vector $w=(w_j)$, obtaining
\begin{align*}
w\T Qw=&\,\frac{\gamma^2T}{\alpha\emph{}}\sum_{k,\ell=1}^m v_m^{[k]}w_k \int_0^1e^{(\lambda_k+
\lambda_\ell)T\tau}d\tau\cdot v_m^{[\ell]}w_\ell\\
=&\,\frac{\gamma^2T}{\alpha}\int_0^1\sum_{k,\ell=1}^m (e^{\lambda_k T\tau}v_m^{[k]}w_k)
( e^{\lambda_\ell T\tau} v_m^{[\ell]}w_\ell) d\tau\\
=&\,\frac{\gamma^2T}{\alpha}\int_0^1\Big(\sum_{k=1}^m e^{\lambda_k T\tau}v_m^{[k]}w_k\Big)^2 d\tau
\ge0.
\end{align*}
This means that $Q$ is semi-definite and $I+Q$ definite and the system \eqref{LGSend} always has a unique solution.
\qed
\par
In general, solutions computed with \eqref{LGSend} will not be sparse, i.e., they will have $m$
nontrivial basis coefficients in the state $y$ and the Lagrange multiplier $p$.
Due the special inhomogeneity in \eqref{RWPv}, sparse solutions for the state $y$ probably do not exist.
However, by adjusting the target profile $\hat y$, one may simply start with a sparse multiplier $p(t)$ with,
for instance, two terms only,
\begin{align}\label{pAnsatz}
 p(t)= \delta_1 e^{\lambda_1(T-t)}v^{[1]}+\delta_2 e^{\lambda_2(T-t)}v^{[2]},
\end{align}
with coefficients $\delta_1,\delta_2$ belonging to some \textit{reasonable} form of the control $u$.
Then, by the boundary condition in \eqref{RWPa}, the corresponding target profile has the form
\begin{align}\label{yTarget}
\hat y=y(T)-\delta_1 v^{[1]}-\delta_2 v^{[2]},
\end{align}
where, by \eqref{LsgOC}, the coefficients of $y(T)$ are given by
\begin{align}\label{coeffYT}
\eta_k(T)=&\,e^{\lambda_k T}\eta_k(0)-\frac{\gamma^2T}{\alpha}v_m^{[k]}
\sum_{\ell=1}^2\gamma_\ell v_m^{[\ell]}\varphi_1\left( (\lambda_k+\lambda_\ell)T\right).
\end{align}
We will use this construction in our numerical example.

\begin{figure}[t]
\centering
\includegraphics[width=6.8cm]{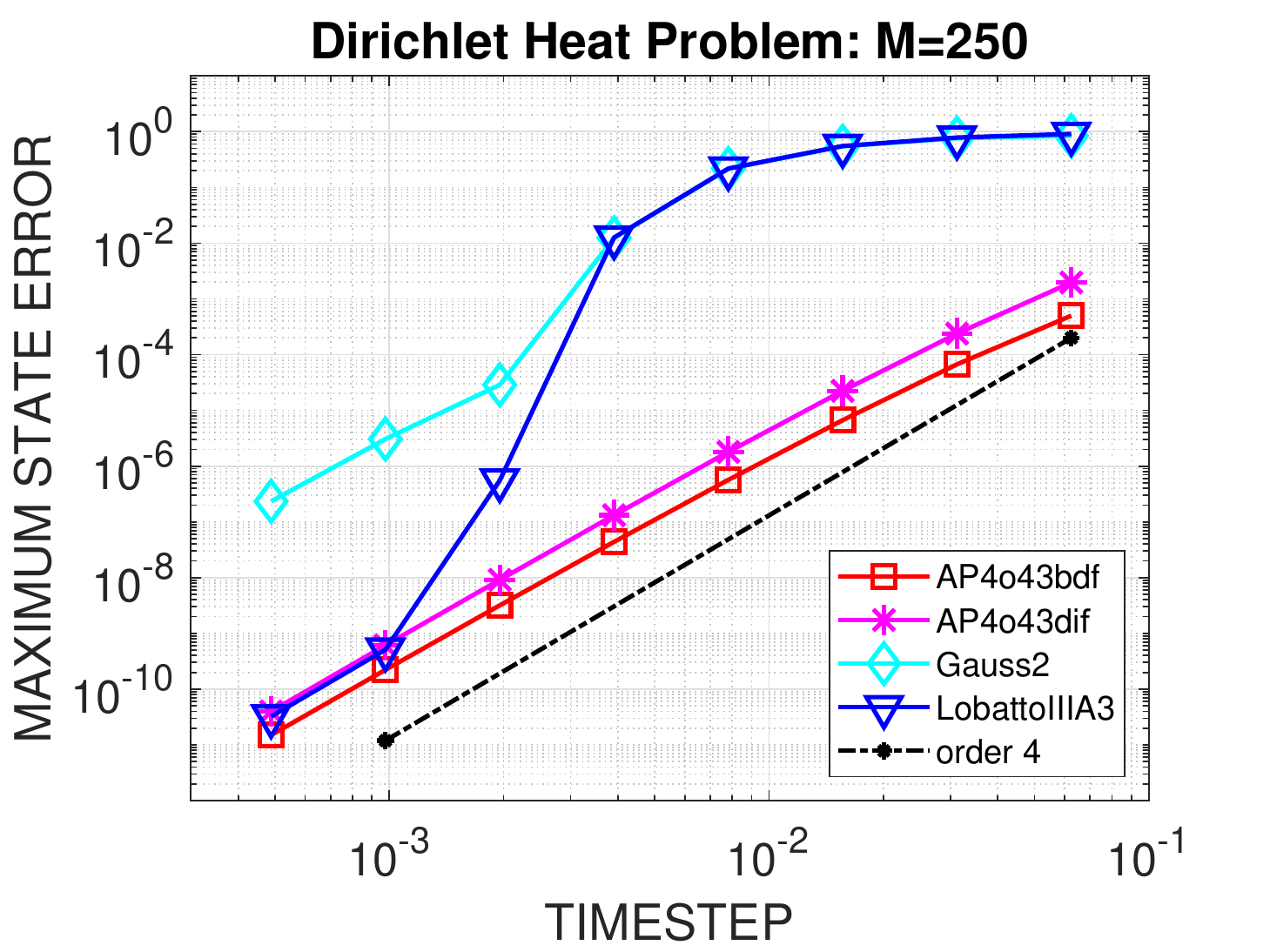}
\hspace{0.1cm}
\includegraphics[width=6.8cm]{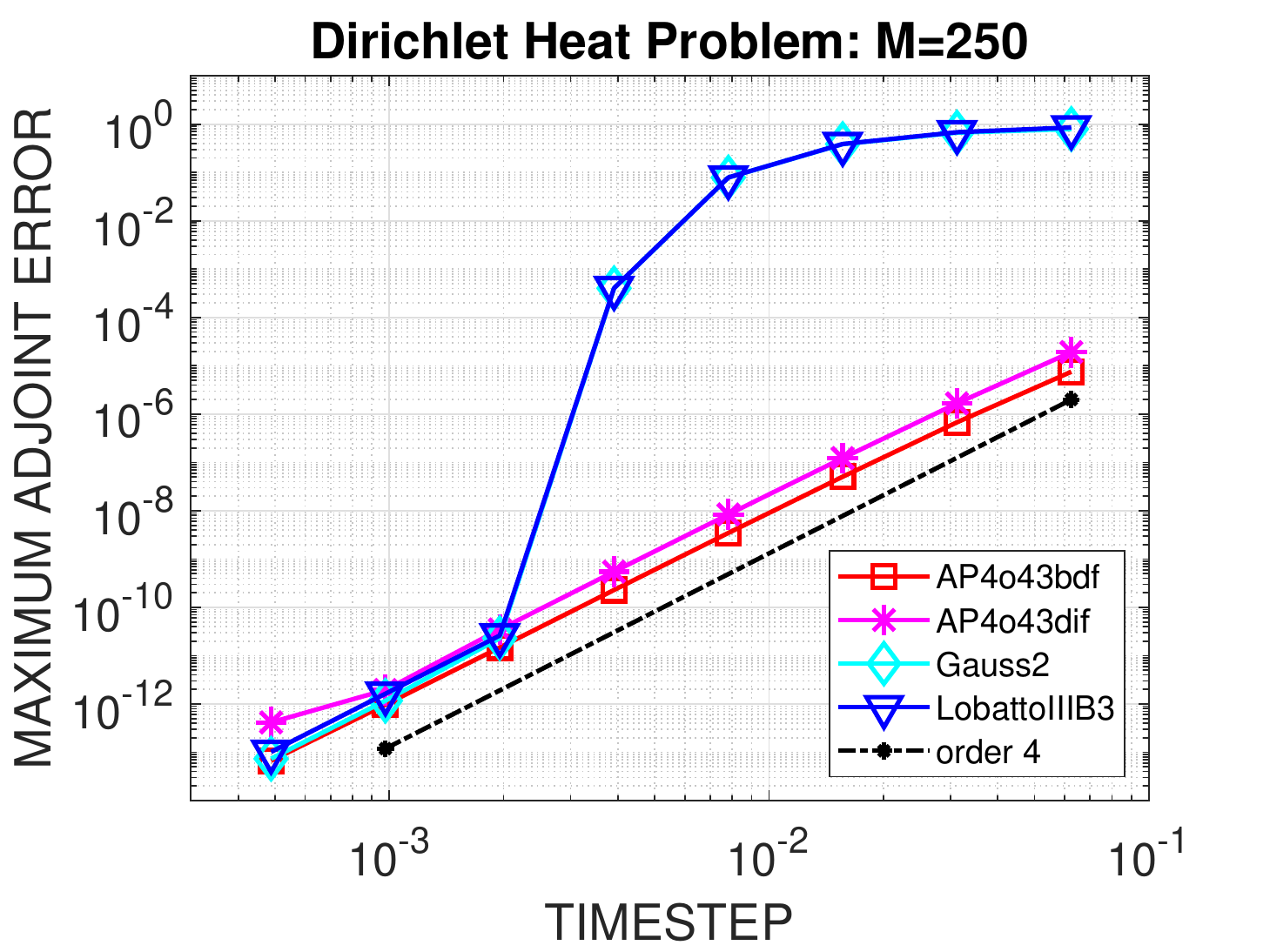}\\[2mm]
\includegraphics[width=6.8cm]{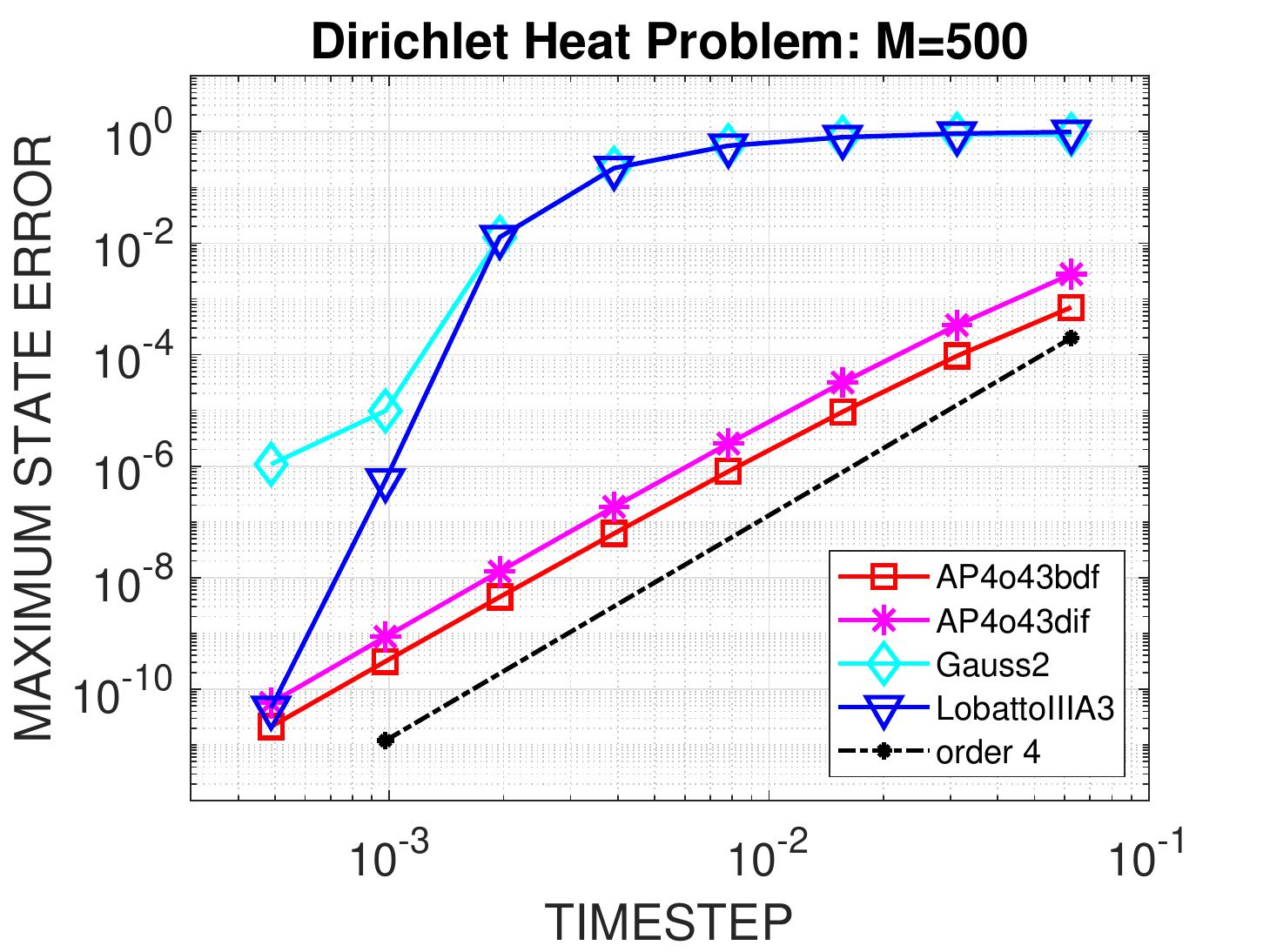}
\hspace{0.1cm}
\includegraphics[width=6.8cm]{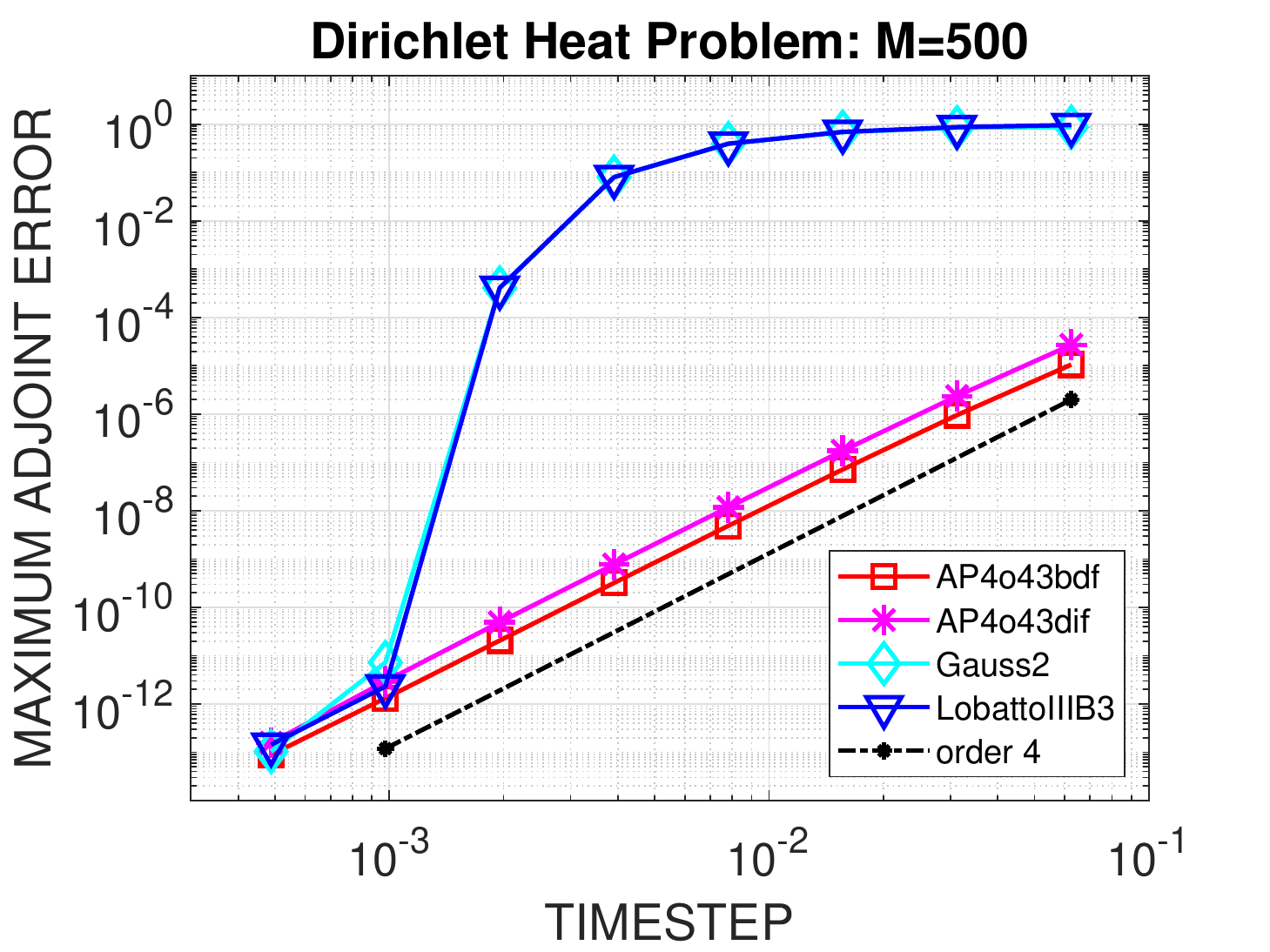}
\parbox{13cm}{
\caption{Dirichlet heat problem with $m=250,500$ spatial points and given exact control $u(t)$:
$\|y(T)$$-$$y_h(T)\|_\infty$ (left) and $\|p(0)$$-$$p_h(0)\|_\infty$ (right).
}\label{fig:yp}
}
\end{figure}

\section{Test case: Dirichlet boundary control problem}
To illustrate an application of the derived expressions for the exact discrete
solutions of the linear heat equation equipped with different boundary conditions, we
consider the following ODE-constrained optimal control problem with an incorporated
boundary control of Dirichlet type:
\begin{align*}
\min_{(y,u)} C&\,:=\frac12\|y(T)-\hat{y}\|^2_2
+\frac{\alpha}{2} \int_0^T u(t)^2\,dt\\
\text{subject to } y'(t)&\, = My(t)+\gamma e_m u(t),\quad t\in (0,T],\\
y(0) &\, = \eins,
\end{align*}
with $T\!=\!1$, $\gamma\!=\!2/\xi^2$, $\alpha\!=\!1$, $\eins=(1,\ldots,1)^T\in\R^m$,
state vector $y(t)\in\R^m$, and $M$ as defined
in \eqref{MatM} with $\theta\!=\!3$. We set $\delta_1\!=\!\delta_2\!=\!1/75$ in
\eqref{pAnsatz} and compute the target profile $\hat{y}\in\R^m$ from \eqref{yTarget}
with coefficients for $y(T)$ defined in \eqref{coeffYT}.
\par
We will compare numerical results for four time integrators of classical order four: the symmetric
2-stage Gauss method \cite[Table II.1.1]{HairerLubichWanner2006}, the symmetric 3-stage partitioned
Runge-Kutta pair Lobatto IIIA-IIIB \cite[Table II.2.2]{HairerLubichWanner2006} and our recently
developed two-step Peer methods AP4o43bdf and AP4o43dif \cite{LangSchmitt2022b}. The two one-step methods
are symplectic and therefore well suited for optimal control \cite{HairerLubichWanner2006,SanzSerna2016}.
\begin{figure}[t]
\centering
\includegraphics[width=6.8cm]{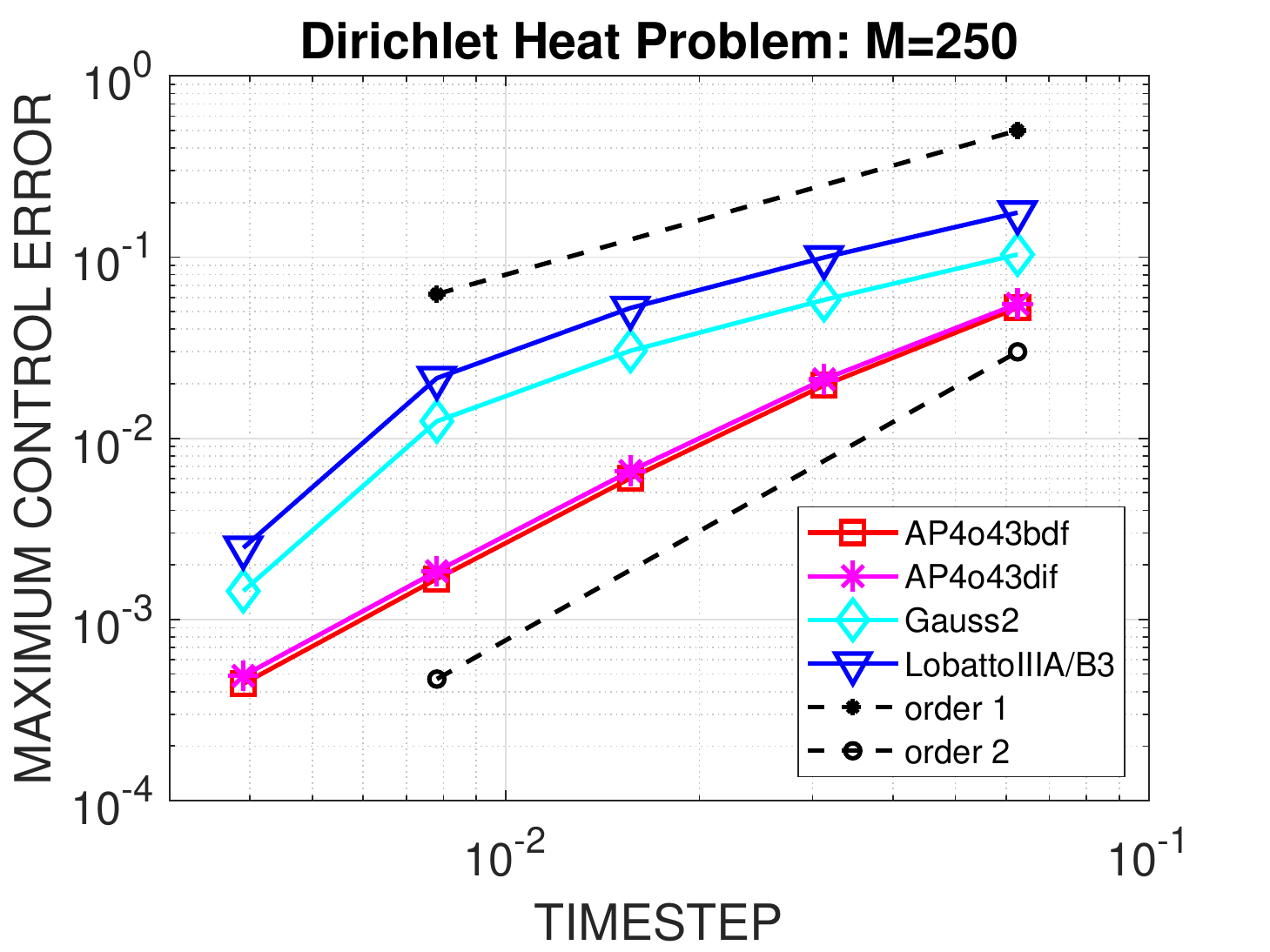}
\hspace{0.1cm}
\includegraphics[width=6.8cm]{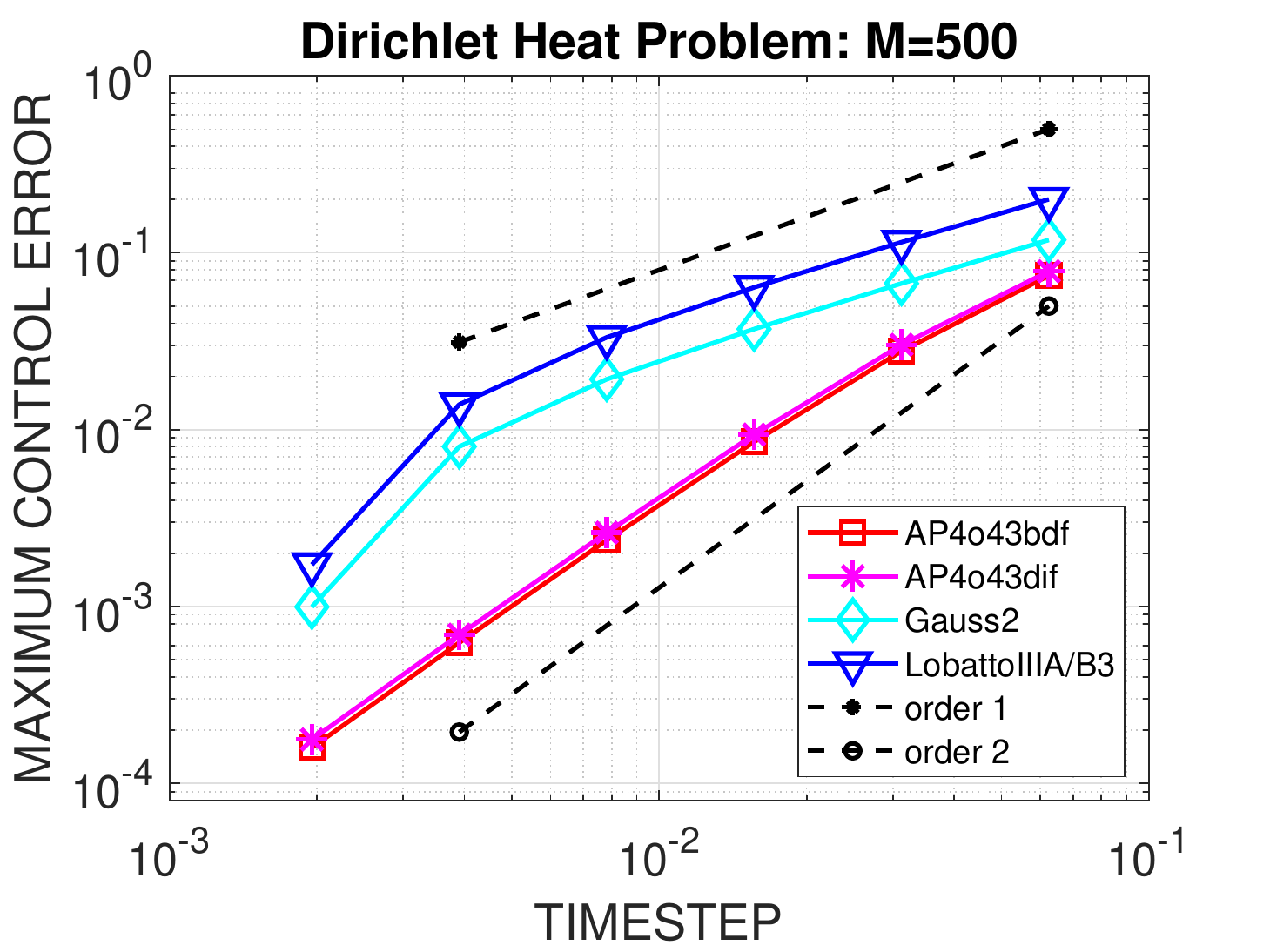}\\[2mm]
\parbox{13cm}{
\caption{Dirichlet heat problem with $m=250,500$ spatial points, solved by gradient descent for $(y,p,u)$: $\max_{n,i}|u(t_{ni})$$-$$u_h(t_{ni})|$
}\label{fig:ctr}
}
\end{figure}
Two test scenarios are considered. First, the accuracy of the numerical approximations for $y(T)$ and
$p(0)$ are studied, where the exact control $u(t)=-\gamma p_m(t)/\alpha$ is used. The initial value for
the multiplier is set to $p(T)=y_h(T)-\hat y$ with $y_h(T)$ being the approximation of $y(T)$ with time step $h$. In this case,
the Karush-Kuhn-Tucker system decouples and only two systems of linear ODEs have to be solved.
In the second scenario, the
optimal control problem is solved for all unknowns $(y,p,u)$ by a gradient method with line search as
implemented in the Matlab routine \textit{fmincon}, see e.g. \cite{Betts2010,Troutman1996} for more details,
and the errors for the control are discussed. The control variable $u(t)$ is approximated at the nodes
$t_{ni}\!=\!t_n+c_ih$, $i=1,\ldots,s,$ chosen by the time integrators on a time grid $\{t_0,\ldots,t_N\}$ with step size $h$ and $s$ stages \cite{Hager2000,LangSchmitt2022b}. In both test cases, we use $m\!=\!250$ and $m\!=\!500$ to also study the influence of the system size. The number of time steps are $N=2^k$ with $k=4,\ldots,11$.
\par
In Fig.~\ref{fig:yp}, results for the first test scenario are shown. Not surprisingly,
the serious order reduction for the symplectic one-step Runge-Kutta methods are clearly seen.
This phenomenon is well understood and occurs particularly drastically for time-dependent
Dirichlet boundary conditions \cite{OstermannRoche1992}. This drawback is shared by all one-step
methods due to their insufficient stage order. Note that the number of affected time steps
increases when the system size is doubled. In contrast, the newly designed two-step Peer methods for optimal control problems
work quite close to their theoretical order four for the state $y$ and the adjoint $p$.
The order reduction for the one-step methods is also visible for the more challenging
fully coupled problem. The results plotted in Fig.~\ref{fig:ctr} show a reduction to first
order for the approximation of the control, whereas the two-step methods perform with
order two for this problem. We refer to \cite{Hager2000} for a discussion of the convergence
order for general ODE constrained optimal control problems. Once again, the range of the affected
time steps depends on the problem size. It increases for finer spatial discretizations.

\bibliographystyle{plain}
\bibliography{bibpeeropt}
\end{document}